\documentclass[s5,11pt]{article}
\usepackage{amsmath}
\newtheorem{theorem}{Theorem}

\newtheorem{lemma}{Lemma} 

  \numberwithin{equation} {lemma}
\author{ Tord Sj\"odin}
\title{ A short proof of The Fundamental Theorem of Algebra }
 \begin{document}  

\maketitle
\begin{abstract} D'Alembert made the first serious attempt to prove the Fundamental Theorem of Algebra (FTA) in 1746. Although it contained gaps, it became a model for several later proofs by C. F. Gauss and others. An elementary proof of (FTA), based on the same idea,  is given in Proofs from THE BOOK \cite{AZ}. We give a shorter and more transparent version of this proof.\end{abstract}
\paragraph{ \it AMS 2010 Subject Classsification:} Primary 12 - 01
Secondary 30 - 01
\paragraph{ \it Key words and phrases:} Complex numbers, polynomial, coefficient, minimum point, root 
\section{Introduction.} The Fundamental Theorem of Algebra (FTA) is one of the truly classical theorems in modern mathematics, with a history that goes back at least to the beginning of the 18-th century. 
It states that every polynomial of degree at least one and with complex coefficients has at least one complex zero. A useful consequence is that a complex polynomial of degree $n\geq 1$ has exactly $n$ complex zeroes, counting multiplicities. This stands in contrast to real polynomials that may lack real zeroes. The reason for this difference between real and complex polynomials is usually expressed by saying that the set of complex numbers, but not the set of real numbers, is algebraically closed.\\[1em]
An elementary proof of (FTA) based on d'Alembert's idea was given by Fefferman \cite{F} and in shorter form by Aigner, Ziegler \cite{AZ}. There are many other more or less elementary proofs of (FTA), see \cite{B}, \cite{R} and \cite{W}.  We improve on the proof in \cite{AZ} by a much shorter and more transparent argument in the key lemma. As in \cite{AZ}, we will take the following three facts as known:\\[0.5em]
(a) {\it Polynomials and their absolute values are continuous,}\\[0.5em]
(b) {\it A continuous function on a compact set attains a minimum value,}  \\[0.5em]
(c) {\it Every complex number with absolute value one has an m-th root}, $m\geq 1$.
\begin {theorem} Every non-constant complex polynomial $p(z)=a_nz^n+a_{n-1}z^{n-1}+\cdots +a_1z+a_0$ has at least one root.\end{theorem} 
The following lemma is used by many authors. Our contribution is a more effective, transparent and geometric construction of the critical set.\begin{lemma} Let $p(z)$ be a non--constant polynomial with $p(a)\ne 0$, then every open disc $D$ around $a$ contains a  point $b$ with $|p(b)|<|p(a)|$.
\end{lemma}
Proof. Let $p(z)$ be as in the lemma. We first claim that, without loss of generality, we may assume that $a=0$ and $p(a)=1$. If this is not the case, we define another polynomial $q(z)=p(z+a)/p(a)$ such that $q(0)=1$, see \cite{R}. Now assume that every open disc $B(0,R)$ of radius $R>0$ around the origin contains a point $b$ such that $|q(b)|<1$, then the open disc $D=B(a, R)$ around $a$ contains the point $a+b$ such that $|p(a+b)|<|p(a)|$, which proves our claim.\\[1em]
In the following we assume that $p(z)=1+a_1z+a_2z^2+\cdots +a_nz^n$ and let $m$ be the smallest integer such that $a_m\ne 0$. Then $p(z)$ can be written in the form
\begin{displaymath}p(z)=1+a_mz^m+ z^{m+1}(a_{m+1}+\cdots +a_nz^{n-m-1})=1+a_mz^m+r(z).\end{displaymath}
The rest of the proof is in two steps. In the first step we find $0<\rho <1$ such that
\begin{equation}|r(z)|<|a_mz^m|<1,\quad\textrm{for all }\,0< |z|\leq \rho .\end{equation}
The second inequality in (1.1) holds if $|z|< |a_m|^{-1/m}=\rho _1$. To get the first inequality in (1.1) we note that if $|z|<1$ then
\begin{displaymath} |r(z)|\leq |z|^{m+1}(|a_{m+1}|+\cdots +|a_n|)<|a_m|\cdot |z|^m=|a_mz^m|
\end{displaymath}
provided $0<|z|<(|a_{m+1}|+\cdots |a_n|)/|a_m|=\rho _2$. We conclude that (1.1) is true for any $0<\rho <\min (\rho _1,\rho _2,1)$ and fix any such $\rho$.\\[1em]
In the second step we let $\zeta$ be an m-th root of the complex number $-\overline{a_m}/|a_m|$ and put $w=\rho\cdot\zeta$. Then $|w|=\rho<1$, 
\begin{displaymath}a_mw^m=-a_m\cdot \rho ^m\cdot \frac{\overline a_m}{|a_m|}=-|a_m|\cdot \rho ^m
\end{displaymath}
by the definition of $w$ and $|r(w)|<|a_m|\cdot \rho ^m<1$, by (1.1). Hence
\begin{displaymath}|p(w)|\leq |1+a_mw^m|+|r(w)|=1-|a_m|\rho ^m+|r(w)|<1.\end{displaymath}
Since $\rho$ can be taken arbitrarily small, every open disc $D$ around the origin contains some point $b$ where $|p(b)|<1$. This completes the proof of the lemma.\hfill $\triangle$\\[1em]
{\it Remark.} The last inequality has a simple geometric interpretation. Let $B(0,1)$ be the unit disc and let $x =1-|a_m|\rho ^m$. Then $x$ is a point between $0$ and $1$ on the real axis and $p(w)$ lies in the open disc $B(x,1-x)\subset B(0,1)$.
\\[1em]
{\it Proof of the theorem}. Let $p(z)=a_nz^n+a_{n+1}z^{n+1}+\cdots +a_1z+a_0$ be an n-th degree polynomial with $a_0\ne 0$, $a_n\ne 0$ and let $\overline B(0,R)$ denote the closed disc with radius $R$ centered at the origin. Then the continuous function $|p(z)|$ attains a smallest value on the compact set $\overline B(0,R)$ at $z_0$. We claim that $z_0$ is an interiour point of $\overline B(0,R)$, provided $R$ is choosen large enough. To see this, we have for $|z|\geq 1$ the estimate
\begin{displaymath}|p(z)|\geq |a_n|\cdot |z|^n-|z|^{n-1}  (|a_{n-1}|+\cdots  +|a_0|)\geq |a_n|\cdot |z|^n/2 \rightarrow \infty ,
\end{displaymath}
as $|z|\rightarrow \infty$. Then we can choose $R$ such that $|p(z)|>|p(0)|\geq |p(z_0)|$, for $|z|=R$, which proves our claim. Now assume that $p(z_0)\ne 0$. Then by the lemma there is an open disc $B(z_0,\rho )\subset B(0,R)$ and a point $b$ in $B(z_0,\rho )$ such that $|p(b)|<|p(z_0)|$. This contradicts the minimum property of $z_0$ and proves that in fact $p(z_0)=0$.\hfill $\triangle$

 \end{document}